\documentclass[12pt]{amsart}

\usepackage{amsfonts}
\usepackage{amssymb}
\usepackage{amsthm}
\usepackage{enumerate}
\usepackage{hyperref}

\def\fracd#1#2{\frac{\displaystyle #1}{\displaystyle #2}}
\def\Sp{{\rm Sp}}

\newcommand\Ad{{\rm Ad}}
\newcommand\ad{{\rm ad}}

\newcommand\C{{\mathbb C}}

\newcommand\Cl{{\rm C \ell}}

\newcommand\End{{\rm End}}
\newcommand\eS{{\sf S}}
\newcommand\Fix{{\rm Fix}}

\newcommand\g{{\mathfrak g}}

\newcommand\HH{{\sf H}}
\newcommand\Hom{{\rm Hom}}
\newcommand\id{{\rm id}}
\newcommand\II{{\rm I\hspace{-.02em}I}}

\newcommand\K{{\mathbb K}}

\newcommand\Lb{{\mathfrak b}}

\newcommand\LG{{\rm G}}

\newcommand\m{{\mathfrak m}}

\newcommand\qH{{\mathbb H}}
\newcommand\qO{{\mathbb O}}
\newcommand\R{{\mathbb R}}

\newcommand\sgn{{\rm sgn}}

\newcommand\so{{\mathfrak s}{\mathfrak o}}
\newcommand\SO{{\rm SO}}
\newcommand\spann{{\rm span}}
\newcommand\spin{{\mathfrak s}{\mathfrak p}{\mathfrak i}{\mathfrak n}}
\newcommand\Spin{{\rm Spin}}

\newcommand\su{{\mathfrak s}{\mathfrak u}}
\newcommand\SU{{\rm SU}}

\newcommand\T{{\rm T}}
\newcommand\U{{\rm U}}

\renewcommand\k{{\mathfrak k}}
\renewcommand\ker{{\rm ker} \, }

\renewcommand\O{{\rm O}}
\renewcommand\P{{\sf P}}

\renewcommand\sp{{\mathfrak s}{\mathfrak p}}

\renewcommand\u{{\mathfrak u}}

\newtheorem{theorem}[equation]{Theorem}
\newtheorem{lemma}[equation]{Lemma}
\theoremstyle{definition}

\theoremstyle{remark}
\newtheorem{remark}[equation]{\bf Remark}
\numberwithin{equation}{section} \theoremstyle{corollary}
\newtheorem{corollary}[equation]{Corollary}
\theoremstyle{position}
\newtheorem{proposition}[equation]{Proposition}

\newcounter{refcounter}
\renewcommand\therefcounter{\bf (\arabic{refcounter})}
\newcommand\mylabel[1]{\refstepcounter{refcounter}\bf\therefcounter\label{#1}}

\begin{document}

%
%
%

\title{Manifolds with large isotropy groups}

\author{Andreas Kollross and Evangelia Samiou}
\thanks{Research supported by the University of Cyprus}
\address{Institut f\"ur Mathematik\\Universit\"at Augsburg\\86135 Augsburg\\Germany}
\email{kollross@math.uni-augsburg.de}

\address{University of Cyprus\\Department of Mathematics and Statistics\\P.O.\ Box 20537\\1678 Nicosia\\Cyprus}
\email{samiou@ucy.ac.cy}

\subjclass[2000]{53C30, 22E46, 57S15}

\date{\today}

\keywords{homogeneous spaces, isotropy representation, low
cohomogeneity}

\begin{abstract}
We classify all simply connected Riemannian manifolds whose isotropy
groups act with cohomogeneity less than or equal to two.
\end{abstract}


\maketitle


    \section{Introduction and results}
    \label{intro}


Let $M$ be a simply connected Riemannian manifold and $G$ a closed
connected subgroup of the isometry group of~$M$. We consider the
action of an {\em isotropy group}, i.e.\ the subgroup $G_p = \left
\{ g \in G \mid g \cdot p = p \right \}$ consisting of all elements
which leave one point $p \in M$ fixed. This group $G_p$ acts on $M$
by restriction of the $G$-action and by the differentials of the
isometries in $G$ on the tangent space $\T_p M$; we call the first
the {\em isotropy action}, while the latter is called {\em isotropy
representation}. For any differentiable Lie group action, the
minimal orbit codimension is called the {\em cohomogeneity} of the
action.
\par
It is a well known classical result that {\em two-point homogeneous
spaces}, i.e.\ Riemannian manifolds where any ordered pair of
equidistant points can be mapped to any other such pair by an
isometry, are exactly the rank-one symmetric and Euclidean
spaces~\cite{helgason}. Two-point homogeneity is equivalent to the
condition that the isotropy action at any point is of cohomogeneity
one.
\par
Thus, Riemannian manifolds whose isotropy actions are of
cohomogeneity one are well-known and it is a natural question to ask
which Riemannian manifolds have isotropy actions of low
cohomogeneity. In this article, we classify all simply connected
Riemannian manifolds~$M$ whose isotropy groups act with
cohomogeneity less than or equal to two. For homogeneous spaces this
condition is equivalent to saying that the unit tangent bundle
of~$M$ is of cohomogeneity one. However, there are also
non-homogeneous examples of such manifolds.
\par
By $K = G_{p}$ we will usually denote the isotropy subgroup of the
$G$-action at a point~$p$ lying in a principal orbit. Furthermore we
denote by~$c(G,M)$ the cohomogeneity of the $G$-action on~$M$. For
the cohomogeneity of the $G$-action on the tangent bundle we have
\begin{equation*}
c(G,TM)=c(K,T_pGp) + 2 c(G,M) = c(K,M) + c(G,M).
\end{equation*}

Our main result for the homogeneous case is given by the following
theorem. For the non-homogeneous case, see Theorem~\ref{InhomTh}. Of
course, if we allow the isotropy group to be trivial, we get all
two-dimensional manifolds.


\begin{theorem}\label{Coh2Result}
Let $M = G / K$ be a simply connected Riemannian homogeneous space
such that the isotropy group~$K$ acts with cohomogeneity two on~$M$.
Assume further that $M$ is de~Rham irreducible. Then $M$ is
isometric to one of the following homogeneous spaces equipped with a
$G$-invariant metric.
\begin{enumerate}

  \item
  An irreducible symmetric space of rank two.

  \item
  A generalized Heisenberg group of type
  $$
    N(1,n),\; N(2,n),\; N(3;n,0),\; N(6,1),\; N(7;1,0).
  $$

  \item
  One of the following compact homogeneous spaces:
  $$
    \frac{\SU(n+1)}{\SU(n)},\;
    \frac{\Sp(n+1)}{\U(1)\cdot\Sp(n)},\;
    \frac{\Sp(1)\cdot\Sp(n+1)}{\Delta\Sp(1)\cdot\Sp(n)},\;
    \frac{\Spin(9)}{\Spin(7)}.
  $$

  \item
  One of the following non-compact homogenous spaces:
  \begin{align*}
    \fracd{\SU(n,1)}{\SU(n)},\;
    \fracd{\Sp(n,1)}{\U(1)\cdot\Sp(n)},\;
    \fracd{\Sp(1)\cdot\Sp(n)\ltimes\R^{4n}}{\U(1)\cdot\Sp(n)},\;
    \fracd{\Sp(1)\cdot\Sp(n,1)}{\Delta\Sp(1)\cdot\Sp(n)},\qquad\;\\
    \fracd{\Sp(1)\cdot(\Sp(1)\cdot\Sp(n)\ltimes \R^{4n})}{\Delta
    \Sp(1)\cdot\Sp(n)},\
    \fracd{\Spin(7)\ltimes \delta_7}{\Spin(6)},\;
    \fracd{\Spin(8,1)}{\Spin(7)},\;
    \fracd{\Spin(8)\ltimes\delta_8^+}{\Spin(7)}.
  \end{align*}

  \item
  Hyperbolic space with $G$-action as given in
  Theorem~\ref{nisomisotr}.

\end{enumerate}
\end{theorem}

The paper is organized as follows. We start with a splitting
criterion for homogeneous spaces which we use to detect Riemannian
products from the isotropy representation. In
Section~\ref{isostruct} we prove fibration theorems for homogeneous
spaces with reducible isotropy representation with two submodules.
The actual classification of homogeneous spaces with isotropy
cohomogeneity two is done in Section~\ref{HomClass}. We first
determine all effective representations of cohomogeneity two. Using
the splitting criterion from Section~\ref{derham}, we eliminate
those leading to Riemannian products. Then we construct all
homogeneous spaces corresponding to the remaining representations.
The final Section~\ref{inhomclass} contains the classification in
the inhomogeneous case.

We would like to thank Carlos Olmos for helpful discussions.

    \section{A splitting criterion for homogeneous spaces}
    \label{derham}


If $M$ is isometric to a product $M=X\times Y$, with $X$ irreducible
or maximal Euclidean  then the factors $X$ and $Y$ are preserved by
any connected group of isometries. In particular, the action of any
connected group of isometries on $M$ projects to actions of lower
cohomogeneity on both factors $X$ and $Y$. Furthermore, if $M$ is
homogeneous then so are $X$ and $Y$. For the classification we may
therefore exclude reducible Riemannian manifolds. These can be
detected by means of Corollary~\ref{split} from the isotropy
representation.


\begin{lemma}\label{tgfs}
Let $M=G/K$ be a Riemannian homogeneous space and $V\subseteq T_pM$
a submodule of the isotropy representation. Assume that $V$ is the
fixed point set of some subgroup $N \subseteq K$. Then the
$G$-invariant distribution $\mathcal{D}_V$ given by
$$\mathcal{D}_V(g p)=d_pgV$$
is well-defined and integrable. The leaf through $p \in M$ is a
connected component of the fixed point set of~$N$ and therefore the
leaves of the foliation corresponding to $\mathcal{D}_V$ are closed
totally geodesic submanifolds.
\end{lemma}


\begin{proof}
It follows from the $K$-invariance of~$V$ that the
distribution~$\mathcal D$ is well-defined. We may assume that $N
\subseteq  K$ is the kernel of the $K$-action $K \to \SO(V)$ on~$V$
and let $F \subseteq M$ be the connected component containing~$p$ of
the fixed point set of the $N$-action on $M$. Thus $F \subset M$ is
a totally geodesic submanifold and $T_p F = V$.

We are going to show that $F$ is everywhere tangent to the
distribution $\mathcal D$. For $gp=q \in F$, let
$$
V_g = g^{-1} T_q F\subseteq T_pM \quad\mbox{and}\quad
N_g=g^{-1}Ng\subseteq K\ .
$$
Clearly, $N_g$ fixes $V_g$ pointwisely, i.e.\ $nu=u$ for each $n \in
N_g$ and $u \in V_g$. Also, since $N_g\subseteq K$, it preserves the
decomposition $T_pM=V\oplus V^\perp$ of~$T_pM$ in $K$-invariant
submodules, i.e.\ $N_g V = V$ and $N_g V^\perp = V^\perp$. We have
constructed a continuous map
$$
g \mapsto V_g \in G_{\dim(V)}(V \oplus V^\perp)
$$
into the Grassmannian sending the identity element $e \in G$ to
$V_e=V$. A chart around $V$ for this Grassmannian is given by
\begin{equation}\label{chartGr}
\Hom(V,V^\perp) \to G_{\dim(V)}(V \oplus V^\perp), \quad A \mapsto
\Gamma(A)= \{(v, Av) \mid v \in V \}.
\end{equation}
Hence for $g$ sufficiently close to $e$ we can write
$$
V_g = \{(v, A_g v) \mid v \in V \}
$$
with $A_g \in \Hom(V,V^\perp)$ depending continuously on~$g$ and
$A_e = 0$. Since $N_g$ fixes $V_g$ pointwisely we have for $n \in
N_g$ and $v \in V$ that
\begin{equation}\label{coeffcompx89i}
(nv, n A_g v) = n (v, A_g v) = (v, A_g v)
\end{equation}
and thus $nv = v$, and $n A_g v = A_g v$. Therefore $N_g \subseteq
N$ and, since $N_g$, $N$ are isomorphic compact Lie groups, we
actually have $N_g = N$. From \eqref{coeffcompx89i} we now have $n
A_g v = A_g v$ for all $n \in N$. Thus $A_g v \in V \cap V^\perp = 0
$.

We have shown that $V_g = V$ as long as $V_g$ is in the image $U_0$
of the map \eqref{chartGr}. Since $F$ is connected, for each $q \in
F$ there is a path $g(t) \in G $, $t \in [0,1] $ such that $g(t) p
\in F$, $g(1) p = q$ and $g(0) = e$. The set $\{t \in [0,1] \mid
V_{g(t)} \in U_0 \} \subset [0,1]$ is open and closed since, by the
above, ``$ V_{g(t)} \in U_0 $'' is equivalent to ``$ V_{g(t)} = V
$''. Hence $V_g=V$ for all $g \in G$ with $gp \in F$. In particular,
$T_q F = g V = \mathcal{D}_V(q)$.
\end{proof}


Lemma~\ref{tgfs} provides a splitting criterion for simply connected
homogeneous spaces with reducible isotropy representation.


\begin{corollary}\label{split}
Let $M=G/K$ be an effectively presented simply connected Riemannian
homogeneous space. Assume that the isotropy representation
\begin{equation}\label{splitisotr56}
T_pM=\m = \m_1 \oplus \m_2
\end{equation}
splits with nontrivial $K$-submodules $\m_1 $, $\m_2 $. Additionally
suppose that both $\m_1$, $\m_2$ are fixed modules of suitable
subgroups of~$K$, i.e.\ there are subgroups $N_i \subset K$ such
that
$$
 \m_i = \Fix (N_i) = \left \{v \in \m \mid n v = v \mbox{ for all }
 n \in N_i \right \}.
$$
Then $M$ admits a nontrivial de~Rham decomposition $M = M_1 \times
M_2$.
\end{corollary}


\begin{proof}
We have subgroups $N_i =\ker \left ( K \to \SO(\m_i) \right )
\subset K$ such that
$$
\m_i = \{x \in \m \mid n x = x \mbox{ for all } n \in N_i\}.
$$
By Lemma~\ref{tgfs} we have complementary integrable totally
geodesic $G$-in\-va\-ri\-ant distributions~$\mathcal{D}_{\m_i} $, $i
= 1,2$,
\begin{equation}\label{splitdistr56}
TM= \mathcal{D}_{\m_1} \oplus \mathcal{D}_{\m_2}.
\end{equation}
By $G$-invariance, these distributions are perpendicular at each
point $q \in M$. The distributions are also parallel. To see this,
let $x,y \in\Gamma\mathcal{D}_{\m_1}$ and $v,w
\in\Gamma\mathcal{D}_{\m_2}$ be sections and compute
\begin{equation*}
\begin{split}
   \langle\nabla_x v\mid y\rangle&= -\langle v\mid\nabla_x y\rangle
 = -\langle v\mid \II_1(x,y)\rangle =0\quad\mbox{ and } \\
   \langle\nabla_w v\mid y\rangle&= \langle \II_2(w,v) \mid y\rangle
 =0
\end{split}
\end{equation*}
where $\II_i$ denotes the second fundamental form of
$\mathcal{D}_{\m_i} $. This shows that $\nabla_{x+w} v \in
\mathcal{D}_{\m_i}$ for each $x+w \in TM$ and each section $v \in
\Gamma\mathcal{D}_{\m_i}$. Thus \eqref{splitisotr56} is preserved by
the holonomy of~$M$ at $p$. By the de~Rham Holonomy Theorem $M = M_1
\times M_2$ splits as a Riemannian product according to
\eqref{splitdistr56}.
\end{proof}


    \section{Splitting isotropy and structure
             of homogeneous spaces}
    \label{isostruct}


\begin{theorem}\label{g1}
Let $M = G / K$ be a Riemannian homogeneous space whose isotropy
representation can be decomposed into two (not necessarily
irreducible) invariant subspaces
$$
T_pM = \m = \m_1 \oplus \m_2
$$
such that the kernel~$N_1$ of the $K$-representation on~$\m_1$ acts
without nonzero fixed vectors on~$\m_2$. Then $\g_1 = \k \oplus
\m_1$ is a subalgebra of the Lie algebra~$\g$ and the corresponding
connected subgroup $G_1$ of~$G$ is closed.
\end{theorem}


\begin{proof}
We first show that $\g_1 = \k \oplus \m_1$ is a Lie algebra. Clearly
$[\k,\k] \subseteq \k \subseteq \g_1$ and $[\k,\m_1] \subseteq \m_1
\subseteq \g_1$. It remains to show that $[\m_1,\m_1] \subseteq
\g_1$, or, equivalently, that for $m, m^\prime \in \m_1$ we have $
[m, m^\prime] \perp \m_2$. To this end, let $\mathrm{pr} \colon \g
\to \m_2$ be the orthogonal projection. Since $N_1 \subseteq K $ and
$K$ acts isometrically on~$\g$, this map is $N_1$-equivariant. Now
for $n \in N_1$ we compute
$$
\Ad_n\, \mathrm{pr}\, [m,m^\prime] = \mathrm{pr}\, \Ad_n\,
[m,m^\prime] = \mathrm{pr}\, [\Ad_n\, m,\Ad_n\, m^\prime] =
\mathrm{pr\,} [m,m^\prime] \ .
$$
Hence $\mathrm{pr} [m,m^\prime] $ is a fixed vector for $N_1$ in
$\m_2$ and therefore zero.

In order to show that the connected subgroup $G_1 \subseteq G$ with
Lie algebra~$\g_1$ is closed, recall that in the isotropy
representation, $\m_1$ is the fixed point set of~$N_1$. Thus by
Lemma~\ref{tgfs} the distribution $\mathcal{D}_{\m_1} $ is
integrable with closed totally geodesic leaves. The orbit~$G_1p $ is
an immersed submanifold of~$M$. We also have that
$\mathcal{D}_{\m_1}(p) = \m_1 = T_p G_1 p$ and, by left invariance
of the distribution $\mathcal{D}_{\m_1}$, that
$\mathcal{D}_{\m_1}(q) = T_qG_1 p$ for each $q = g_1 p \in G_1 p$.
Since $G_1p $ is connected, it is contained in the leaf $F$ of
$\mathcal{D}_{\m_1} $ through $p$. We now show that $G_1p = F$.
Since both have the same dimension $\dim \m_1$, there is $r > 0$
such that $F \cap B_r(p) \subseteq G_1p$, where $B_r(p) = \left \{ x
\in M \mid d(x,p) <r \right \}$. If $q = g_1 p \in G_1 p$, then $g_1
F = F$ and $F \cap B_r(q) = g_1(F \cap B_r(p)) \subseteq g_1 G_1 p=
G_1 p $. Thus $G_1 p$ has the following property:
$$
q \in G_1 p \Rightarrow F \cap B_r(q) \subseteq G_1p
$$
But then also the complement $F \setminus G_1 p$ has this property.
It follows that both $G_1p $ and $F \setminus G_1 p$ are open
in~$F$. Since $F$ is connected, we have $G_1 p=F $. Now $G_1$ is
closed, since it is the preimage $G_1 = \mathrm{pr}^{-1} (G_1 p)$
under the projection $\mathrm{pr} \colon G \to G / K = M $,
$\mathrm{pr}\colon g \mapsto g p = g K \in M $.
\end{proof}


In general, the adjoint representation of~$G$ restricted to~$G_1$
does not preserve~$\m_2$, but induces an action on~$\g/\g_1$ by
$$G_1\times\g/\g_1\ni (x,m+\g_1)\mapsto\Ad_x m +\g_1 \in\g/\g_1 \ .$$
Under the isomorphism
\begin{equation*}\begin{split}
\m_2    &\cong\g/\g_1   \\
m   &\mapsto m+\g_1
\end{split} \end{equation*}
and the orthogonal projection $\mathrm{pr} \colon \g \to \m_2$ this
induces a $G_1$-action on~$\m_2$ given by
\begin{equation}\label{projg1act}
G_1\times\m_2\ni(x,m) \mapsto \mathrm{pr} (\Ad_x m)
\end{equation}
which extends the action $K \to \SO(\m_2)$ of~$K$. Henceforth we
will refer to this action simply as the $G_1$-action on~$\m_2$.


\begin{theorem}\label{isomisotr}
Let $M = G / K$ be as in Theorem~\ref{g1} and assume the
$G_1$-action on~$\m_2$ is isometric. Then the quotient $B = G / G_1$
is a Riemannian homogeneous space and we have a Riemannian
submersion and homogeneous fibre bundle
\begin{equation*}
F = G_1 / K \hookrightarrow M = G / K
\stackrel{\pi}{\longrightarrow} G / G_1=B
\end{equation*}
with structure group $G_1$, base $B$ and totally geodesic fibres
$F$.
\end{theorem}


\begin{proof}
Assume that $G_1$ acts isometrically on~$\m_2$. The scalar product
on $\m_2\cong\g/\g_1\cong T_{G_{1}}B$ then extends to a left
invariant metric on $B=G/G_1 $ and the projection $ G/K\to G/G_1$
becomes a Riemannian submersion.
\end{proof}


Let us now look at the case where $G_1$ does {\em not} act
isometrically on~$\m_2$.


\begin{theorem}\label{nisomisotr}
Let $M = G / K$ be as in Theorem~\ref{g1} and such that the
$G_1$-action on~$\m_2$ is {\em not} isometric. Assume that $K$ is
irreducible on~$\m_2$ and for all $m \in \m_2$ the isotropy subgroup
$K_m$ of the $K$-action on~$\m_2$ acts irreducibly on~$\m_1$. Then
$\m_1$ is one-dimensional and $M$ is a semidirect product
\begin{equation}\label{semdirmetr}
M=\R\ltimes\K^n \, \quad\mbox{where}\quad \K=\R,\, \C,\, \mbox{ or
}\qH
\end{equation}
endowed with a left invariant metric arising from the orthogonal
direct sum of standard scalar products on $\R$ and $\K^n$. In
\eqref{semdirmetr} the $\R$-factor acts on~$\K^n$ via a group
monomorphism $\R\to\K^\times$, $t\mapsto e^{\lambda t}$, for some
not purely imaginary $\lambda \in\K$ and $\K^n$ carries the standard
scalar product.
\smallskip\par\noindent
In particular, up to rescaling, $M$ is isometric to hyperbolic
space.
\end{theorem}


\begin{proof}
Suppose $G_1$ does not act isometrically on~$\m_2$. Since $G_1$ is
connected, there is $\xi \in\g_1$ such that the endomorphism
of~$\m_2$ mapping $m \mapsto \mathrm{pr}[\xi,m]$ is not
skew-symmetric, i.e.\ such that $\langle m\mid [\xi,m]\rangle \neq
0$ for some, and therefore almost all, $m \in\m_2$ and $\xi
\in\m_1$. Since $\g_1=\k\oplus\m_1$ and $\ad_\k$ is skew-symmetric
on~$\m_2$ we may assume $\xi \in\m_1$. We fix such an $m \in\m_2$,
$|m|=1$ and consider the homomorphism
\begin{equation}\label{is49}
\m_1\to\R\ , \quad \xi\mapsto\langle m\mid [\xi,m]\rangle
\end{equation}
which is not zero. By the assumption, the isotropy group $K_m$ acts
irreducibly on~$\m_1$, and, since \eqref{is49} is $K_m$-equivariant,
Schur's Lemma yields an isomorphism $\m_1 \cong \R$.

We have $G_1 = K \cdot Q$ with a one-dimensional group $Q \cong
\eS^1$ or $\cong  \R_+$. For $x \in Q$ consider the adjoint $\Ad_x$
of~$x$ acting on $\g=\g_1\oplus\m_2 \ni y+m$. We have
\begin{equation}\label{adaut08}
\Ad_x (y+m) = y + A_x m + L_x m
\end{equation}
with homomorphisms $A_x \in\Hom(\m_2,\g_1)$ and $L_x
\in\Hom(\m_2,\m_2)$.  Since $\Ad_x$ commutes with all $\Ad_k$, $k
\in K$, the homomorphisms $A_x$ and $L_x$ are $K$-module
homomorphisms. It follows from Schur's Lemma that
\begin{equation*}
L_x \in \End (\m_2)^K \cong \K = \R , \ \C , \ \qH
\end{equation*}
and that $A_x \colon \m_2 \to \g_1 = \k \oplus \R$ is zero, an
isomorphism $\m_2\cong\R $ or an isomorphism of~$\m_2$ with some
nontrivial irreducible submodule in the adjoint representation
of~$K$. If $\m_2\cong\R$, then $K$ is trivial and $M=G$ is a
two-dimensional Lie group. In this case it follows that $M$ is
isomorphic to $\R_+ \ltimes \R$, where the left factor acts by
multiplication on the right, since there is only one isomorphism
class of non-abelian two-dimensional real Lie algebras.

Thus we may assume $A_x\colon\m_2\to\k$. We have a one-parameter
group
\begin{equation}\label{Qaction}
L \colon Q \to \K^\times\ ,\quad x \mapsto L_x
\end{equation}
into the multiplicative group of the (skew) field $\K $. Let $\xi
\in \m_1 \cong \R $ be an infinitesimal generator of~$Q$.
Differentiating~\eqref{adaut08} we have
$$ \ad_\xi (y+m) = dA(\xi) m + dL(\xi) m $$
with $dL$ nontrivial, since $G_1$ and therefore $Q$ is not isometric
on~$\m_2$. Thus $0 \neq \lambda = dL(\xi) \in \K $ and we can
diagonalize $\ad_\xi$ on~$\g$. In fact, the eigenspaces are
\begin{equation}\label{eigensp334}
\begin{split}
E(\ad_\xi, 0) &=  \g_1=\k\oplus\R,  \\
E(\ad_\xi, \lambda) &=  \{\lambda^{-1} dA(\xi) m + m \mid m \in\m_2
\}\subseteq \k \oplus \m_2.
\end{split}
\end{equation}
Since $Q$ commutes with $K$ we have $\Ad_K \xi =\xi$ and therefore
the eigenspaces \eqref{eigensp334} are $\Ad_K$-invariant. In
particular, we may replace $\m_1\oplus\m_2$ by $\m_1\oplus
E(\ad_\xi, \lambda) $ as an $\Ad_K$-invariant complement of~$\k$ in
$\g$ to describe the metric structure of~$G/K$. In other words, we
may assume w.l.o.g.\ $\m_2 = E(\ad_\xi, \lambda) $, $dA(\xi)=0 $ and
thus $A_x=0$ for all $x \in Q$.

As $G_1$ is not isometric on $\m_2$ the image of the $1$-parameter
group $L$ is not contained in the unit sphere of $\K $. Hence
$Q=\R$, and we find $x \in Q$ such that $L_x=\lambda \in\R\subset\K
$, $|\lambda |\neq 1 $. For this $x $ we consider the automorphism
$\Ad_x$ of the Lie algebra $\g$. For $m_1,m_2 \in\m_2 $ we compute
\begin{equation*}\begin{split}
\Ad_x([m_1,m_2])=& [m_1,m_2]_{\g_1} + \lambda [m_1,m_2]_{\m_2}  \\
= [\Ad_x m_1,\Ad_x m_2] =& \lambda^2 [m_1,m_2]_{\g_1} + \lambda^2
[m_1,m_2]_{\m_2}
\end{split}\end{equation*}
Since $|\lambda| \neq 1$ we conclude that $[m_1,m_2]=0 $ and $\m_2$
is an abelian ideal in $\g$. It follows that $\m_1\oplus\m_2$ is a
solvable ideal in $\g$.

The connected abelian subgroup $M_2 \subset G$ corresponding
to~$\m_2$ cannot have compact factors since its automorphism group
contains~$K$, which acts irreducibly on the Lie algebra~$\m_2$.
Hence $M_2 = \m_2 \cong \K^{\ell} $ is a $\K$-vector space and the
subgroup of~$G$ corresponding to $\m_1 \oplus \m_2$ is the
semidirect product $Q \ltimes_L M_2$, which acts transitively on the
simply connected manifold~$M$ with discrete isotropy groups. It
follows that $Q \ltimes_L M_2$ acts simply transitively on~$M$ and
we may identify~$M$ with the group $Q \ltimes_L M_2$. Hence we have
$M = Q \ltimes_L M_2 = \R \ltimes_L \K^{\ell}$, where $L$ is the
group homomorphism~(\ref{Qaction}).

To prove that $M$ is isometric to hyperbolic space we only need that
$L(t)=e^{\lambda t}\omega(t)$ with $\lambda \in\R$ and a
one-parameter group $\omega(t) \in\O(\m_2)$ of isometries of $\m_2$.
Identifying $M$ with the product $\R\times\m_2$, the group structure
is given by
$$(t,x) \cdot (s,y) = (t+s, x +e^{\lambda t}\omega(t) y) $$
and the left invariant metric is
$$g(t,x)= dt^2 + e^{-2\lambda t} dx^2$$
where $dt^2$ and $dx^2$ are the metrics on $\R$ and $\m_2$ induced
by the scalar products. In particular the metric on $M$ is a warped
product and does not depend on $\omega$. The curvature of $g$ is
readily computed to be $-\lambda^2$.

To see this, we choose vector fields $T=\frac{\partial}{\partial t}$
tangential to the $\R$-factor and $X,Y$ along the $\m_2$-factor,
which are parallel with respect to the product metric $g_0=dt^2 +
dx^2$. We may also assume that $X$ and $Y$ are perpendicular with
respect to both metrics simultaneously. We abbreviate
$f(t)=e^{-\lambda t}$ and write $|X|$ for the norm of $X$ with
respect to the product metric $g_0$. We have
$$g(T,T)=1,\quad g(X,T)=0=g(X,Y),\quad g(X,X)=f^2|X|^2 \ . $$
The various covariant derivatives are
$$
 \nabla_T T, = 0 \quad\nabla_T X = \frac{f'}{f} X, \quad\nabla_X X
  = -f'f |X|^2 T, \quad\nabla_X Y = 0  \ .
$$
Hence, since the vector fields $T$, $X$, $Y$ commute,
\begin{equation}\begin{split}
K(X+T,Y) &= g(R_{X+T,Y} Y , X+T) = g(\nabla_{X+T}\nabla_{Y}Y-\nabla_{Y}\nabla_{X+T}Y,X+T)   \\
&= g(-|Y|^2\nabla_{X+T}f'fT-\frac{f'}{f}\nabla_{Y}Y,X+T)            \\
&= g(-|Y|^2\frac{f'}{f}f'f X - |Y|^2f''f T - |Y|^2 f'^2 T+\frac{f'}{f}f'f |Y|^2 T,X+T)      \\
&= -|Y|^2 |X|^2 f'^2f^2- |Y|^2f''f  - |Y|^2 f'^2 + f'^2 |Y|^2
\end{split}\end{equation}
For $f(t)=e^{-\lambda t}$ this becomes
$$-\lambda^2(|Y|^2 |X|^2f^4+ |Y|^2f^2) =- \lambda^2 |(X+T) \wedge Y|^2 $$
hence the sectional curvature of the plane spanned by $X+T$ and $Y$
is $-\lambda^2$. Since the set of planes of this type is dense, the
sectional curvature of $M$ is constant $-\lambda^2$.
\end{proof}


    \section{The classification in the homogeneous case}
    \label{HomClass}


In this section, we will classify all simply connected Riemannian
homogeneous spaces whose isotropy action is of cohomogeneity two.
Let $M$ be a simply connected Riemannian homogeneous space and $G$
the connected component of the isometry group of~$M$. We denote by
$K = G_p$ be the isotropy subgroup at a point $p \in M$. Then $K$ is
a compact subgroup of~$G$ and since $M = G / K$ is simply connected
it follows that $K$ is connected. Let $\m = T_pM$ be the tangent
space of~$M$ at the point~$p$ with the scalar product $\langle \
\cdot \mid \cdot \ \rangle$ induced by the Riemannian metric of~$M$.
The action of~$K$ on~$\m$ by $\Ad_G |_K$ is orthogonal with respect
to this scalar product.

We assume that the isotropy representation of~$K$ on~$\m$ is of
cohomogeneity two, thus it is either irreducible or has two summands
$\m = \m_1 \oplus \m_2$ both of cohomogeneity one. We start with the
case where $\m$ is irreducible.


    \subsection{Isotropy irreducible spaces}
    \label{isoirr}


Assume $M = G / K$ is isotropy irreducible, i.e.\ $K$ acts
irreducibly on the tangent space~$T_pM$. In case $M$ is non-compact,
the space is symmetric by~\cite{wolf}. Now assume that $M$ is
compact and has isotropy representation of cohomogeneity~two. By a
result in \cite{wolf}, such a compact $M$ is symmetric or $G$ is
simple and compact. Since a cohomogeneity two representation is
polar, it follows from \cite{kp} that $M$ is symmetric and we have
proved:


\begin{theorem}\label{IsoIrrTh}
An isotropy irreducible simply connected homogenous Riemannian
manifold $M = G / K$ with isotropy representation of cohomogeneity
two is a symmetric space.
\end{theorem}


    \subsection{Reducible isotropy representation}
    \label{ct2}


Now assume $M = G / K$ is a simply connected Riemannian homogeneous
space whose isotropy representation
$$
T_pM = \m = \m_1 \oplus \m_2
$$
consists of two irreducible cohomogeneity one summands $\m_1$ and
$\m_2$, such that $K$ acts transitively on $ \eS\m_1\times \eS\m_2$,
where $\eS\m_i$ denotes the unit sphere in~$\m_i$.

Let $N_i$ be the kernel of~$K\to\SO(\m_i)$. The fixed modules
of~$N_1$ on~$\m_2$ and of~$N_2$ on~$\m_1$ are $K$-invariant
subspaces and therefore trivial by effectivity. If both $N_1$ and
$N_2$ are nontrivial, then $M$ splits by Corollary~\ref{split}. But
we may assume that $M$ is de~Rham irreducible, since otherwise $M$
is just the Riemannian product of two manifolds both of which are
two-point homogeneous. We conclude that $K$ acts effectively on at
least one of~$\m_1$ or $\m_2$, say on~$\m_2$.

We will now determine all reducible compact Lie group
representations of cohomogeneity two without trivial submodules and
where at least one submodule is effective.


\begin{proposition}\label{RepClass}

Let $K$ be a compact Lie group and let $\m = \m_1 \oplus \m_2$ be an
orthogonal $K$-representation such that:

\begin{enumerate}

  \item $\m_1$ and $\m_2$ are invariant subspaces,

  \item $K$ acts nontrivially on~$\m_1$;

  \item $K$ acts effectively on~$\m_2$;

  \item the action of~$K$ on~$\m$ is of cohomogeneity two.

\end{enumerate}
Then $K$, $\m_1$ and $\m_2$ are as given in Table~\ref{Tcoh2}.
\def\fracd#1#2{\frac{\displaystyle #1}{\displaystyle #2}}
\begin{table}[h]\rm
\begin{tabular}{|c|c|c|c|c|}\hline
No. & $K$ & $\dim \m_1$ & $\dim \m_2$ & Range \\
\hline\hline
\mylabel{r203} &  $\U(n)$ & $2$ & $2n$ & $n \ge 3$\\
\hline
\mylabel{r204} &  $\U(1)\cdot\Sp(n)$ & $2$ & $4n$ & $n \ge 1$ \\
\hline
\mylabel{r205} &  $\Sp(1)\cdot\Sp(n)$ & $3$ & $4n$ & $n \ge 1$ \\
\hline
\mylabel{r201} &  $\Spin(6) = \SU(4)$ & $6$ & $8$ & \\
\hline
\mylabel{r202} &  $\Spin(7)$ & $7$ & $8$ & \\
\hline
\end{tabular}\vspace{1em}
\caption{Reducible representations of cohomogeneity
two}\label{Tcoh2}
\end{table}
\end{proposition}

\begin{remark}
While for the entries \ref{r205}, \ref{r201} and \ref{r202} of
Table~\ref{Tcoh2} the representations are determined up to
equivalence by the dimensions of the irreducible modules~$\m_1$ and
$\m_2$, the table entries \ref{r203} and \ref{r204} stand for
infinitely many equivalence classes of representations, depending on
the action of the abelian factor on the module~$\m_1$.  In case of
representation~\ref{r203}, the action of $g \in \U(n)$ on $\m_1 =
\C$ is given by multiplication with a power of the determinant~$\det
g^k$, where we may assume that $k$~is a positive integer, since the
representations of~$K$ by $\det g^k$ and $\det g^{-k}$ on~$\m_1$
(considered to be a real vector space) are equivalent. Similarly, in
case of representation~\ref{r204}, the action of $\U(1) \cdot \Sp(n)
= \left \{ zA \mid z \in \C,\, |z|=1,\, A \in \Sp(n) \right \}$
on~$\m_1 = \C$ is given by multiplication with $z^{2k}$, where
$k$~is a positive integer.

We may impose the condition $n \ge 3$ on representation~\ref{r203}
since for $n = 1$ the cohomogeneity is~$3$ and for $n = 2$ the
action is equivalent to~\ref{r204}.
\end{remark}


\begin{proof}
To prove the proposition, one may use the well known classification
of Lie groups acting transitively and effectively on spheres by
isometries; see e.g.\ \cite{besse},~7.13. For the convenience of the
reader, we reproduce it here in Table~\ref{Tcoh1}.
\begin{table}[h]
\begin{tabular}{|c|c|c|}
\hline  Group $K$ & Isotropy $K_q$& $\dim \m_2$  \\
\hline \hline $\SO(n)$ & $\SO(n-1)$ & $n$     \\ \hline $\SU(n)$ &
$\SU(n-1)$ & $2n$    \\ \hline $\Sp(n)$ & $\Sp(n-1)$ & $4n$    \\
\hline $\U(n)$ & $\U(n-1)$ & $2n$      \\ \hline $\Sp(n)\cdot\Sp(1)$
& $\Sp(n-1)\cdot\Sp(1)$ & $4n$ \\ \hline $\Sp(n)\cdot\U(1)$ &
$\Sp(n-1)\cdot\U(1)$ & $4n$ \\ \hline $\LG_2$ & $\SU(3)$ & $7$
\\ \hline $\Spin(7)$ & $\LG_2$ & $8$      \\ \hline $\Spin(9)$ &
$\Spin(7)$ & $16$  \\ \hline
\end{tabular}\vspace{1em}
    \caption{Representations of cohomogeneity one, $n>1$}\label{Tcoh1}
\end{table}
It is straightforward to find for any such representation~$\m_2$ of
a compact Lie group~$K$ all other (not necessarily effective)
$K$-representations~$\m_1$ of cohomogeneity one such that the action
of~$K$ on $\m_1 \oplus \m_2$ is of cohomogeneity two, i.e.\ the
action of~$K$ on~$\m_1$ restricted to a principal isotropy group of
the $K$-action on~$\m_2$ still acts with cohomogeneity one.

Alternatively, one may start with the classification of reducible
cohomogeneity two representations given in \cite{straume} or
\cite{bergmann} and then exclude those representations where both
modules are non-effective.
\end{proof}


Let $\g$ be the Lie algebra of~$G$ and let $\k$ be the subalgebra
corresponding to~$K$. Then $\g$ is of the form
$$
\g = \k \oplus \m_1 \oplus \m_2,
$$
where $\m = \m_1 \oplus \m_2$ is an $\Ad_G|_K$-invariant complement
of~$\k$ in~$\g$ such that the action of~$K$ on~$\m$ is of
cohomogeneity two, $\m_1$ and $\m_2$ are invariant subspaces
and~$\m_2$ is an effective $K$-representation.

In the following we will determine all possible Lie algebra
structures on~$\g$. Since we have already fixed the Lie bracket
on~$\k \times \k$ and on~$\k \times \m$, it remains to determine the
set of all skew symmetric maps $[\cdot,\cdot] \colon \m \times \m
\to \g$ such that the bracket thus defined on $\g \times \g$
satisfies the Jacobi identity.

One major restriction for the possible Lie algebra structure comes
from representation theory, namely from the fact that any Lie
bracket on~$\g$ defines a $K$-equivariant map $\Lambda^2 \g \to \g$,
since the elements of~$K$ act on~$\g$ as Lie algebra automorphisms.

We distinguish three cases: The case where $\m_1$ is a trivial
one-dimensional submodule is treated in Section~\ref{1triv1irr}. The
case where $\m$ is equivalent to representation~\ref{r203} from
Table~\ref{Tcoh2} is studied in Section~\ref{aconstr}. We begin with
the remaining four representations \ref{r204}, \ref{r205},
\ref{r201} and \ref{r202} from Table~\ref{Tcoh2}.


    \subsection{Heisenberg isotropy}
    \label{Heisenberg}


In this section, we give a classification of all real Lie algebras
$\g$ of simply connected Lie groups~$G$ which have a compact
subgroup~$K$ such that the isotropy representation of the
homogeneous space~$G / K$ is of cohomogeneity two, without a trivial
submodule and equivalent to the isotropy representation of a
generalized Heisenberg group endowed with a left invariant metric.
Some generalized Heisenberg algebras occur as special cases of this
construction.


\begin{proposition}
Let $M = G / K$ be a simply connected Riemannian homogeneous space
which is not de~Rham decomposable and such that the action of~$K$ on
the tangent space of~$M$ is equivalent to one of the representations
\ref{r204}, \ref{r205}, \ref{r201} or \ref{r202}. Then $M$ is
isometric to one of the homogeneous spaces as given in
Table~\ref{Tcoh2Spaces} equipped with a $G$-invariant metric.

\begin{table}[h]\rm
\begin{tabular}{|c|c|c|c|c|c|c|}\hline
No.&$K$&$\dim\m_1$&$\dim\m_2$&$F$&$B$&$M$\\
\hline\hline
\ref{r204}&$\U(1)\cdot\Sp(n)$&$2$&$4n$&$\R^2$&$\R^{4n}$&$N(2,n)$\\
&$=\Spin(2)\cdot\Sp(n)$&&&$\eS^2$&$\qH\P^n$&$\fracd{\Sp(n+1)}{\U(1)\cdot\Sp(n)}$\\
&&&&$\eS^2$&$\R^{4n}$&$\fracd{\Sp(1)\cdot\Sp(n)\ltimes\R^{4n}}{\U(1)\cdot\Sp(n)}$\\
&&&&$\eS^2$&$\qH\HH^{n}$&$\fracd{\Sp(n,1)}{\U(1)\cdot\Sp(n)}$
\\\hline
\ref{r205}&$\Sp(1)\cdot\Sp(n)$&$3$&$4n$&$\R^3$&$\R^{4n}$&$N(3;n,0)$\\
&$=\Spin(3)\cdot\Sp(n)$&&&$\eS^3$&$\qH\P^n$&
$\fracd{\Sp(1)\cdot\Sp(n+1)}{\Delta\Sp(1)\cdot\Sp(n}$\\
&&&&$\eS^3$&$\R^{4n}$&$\fracd{\Sp(1)\cdot(\Sp(1)\cdot\Sp(n)\ltimes\R^{4n})}{\Delta\Sp(1)\cdot\Sp(n)}$\\
&&&&$\eS^3$&$\qH\HH^{n}$&$\fracd{\Sp(1)\cdot\Sp(n,1)}{\Delta\Sp(1)\cdot\Sp(n)}$\\
\hline
\ref{r201}&$\SU(4)=\Spin(6)$&$6$&$8$&$\R^6$&$\R^8$&$N(6,1)$\\
&&&&$\eS^6$&$\R^8$&$\fracd{\Spin(7)\ltimes\delta_7}{\Spin(6)}$
\\\hline
\ref{r202}&$\Spin(7)$&$7$&$8$&$\R^7$&$\R^8$&$N(7;1,0)$\\
&&&&$\eS^7$&$\eS^8$&$\fracd{\Spin(9)}{\Spin(7)}$\\
&&&&$\eS^7$&$\R^8$&$\fracd{\Spin(8)\ltimes\delta_8^+}{\Spin(7)}$\\
&&&&$\eS^7$&$\HH^8$&$\fracd{\Spin(8,1)}{\Spin(7)}$\\
\hline
\end{tabular}\vspace{1em}
\caption{Non-symmetric homogeneous spaces whose isotropy
representation is of cohomogeneity two without a trivial
submodule}\label{Tcoh2Spaces}
\end{table}
\end{proposition}


\begin{proof}
Let $\g = \k \oplus \m_1 \oplus \m_2$ be a real Lie algebra and $G$
be the corresponding simply connected Lie group such that the action
of the connected subgroup~$K$ corresponding to the Lie algebra~$\k$
is equivalent to one of the representations \ref{r204}, \ref{r205},
\ref{r201}, or \ref{r202} given in Table~\ref{Tcoh2}. Assume the
homogenous space $M = G / K$ is endowed with a $G$-invariant
Riemannian metric.

In each case the group $K$ is of the form $K = K_0 \cdot K_1$, where
$K_0 = \Spin(n)$, $n = 2$, $3$, $6$, or $7$, respectively; the
$K$-module $\m_1 = \R^n$ is equivalent to the standard
representation of $\Spin(n)$. Since there is an element in the
center of~$K$ which acts trivially on~$\m_1$ and as minus identity
on~$\m_2$, we know that $\g_1 = \k \oplus \m_1$ is a {\em symmetric
subalgebra} of~$\g$, i.e.\ a (not necessarily compact) subalgebra
of~$\g$ for which the Cartan relations hold; i.e.\ we have
\begin{equation}\label{CartanRel}
[\g_1, \g_1] \subseteq \g_1,\quad [\g_1, \m_2] \subseteq \m_2,\quad
[\m_2, \m_2] \subseteq \g_1.
\end{equation}
We start by determining all possible brackets on~$\m_1 \times \m_1$.
Since the action of~$K$ on~$\m_1$ is the standard representation
of~$\Spin(n)$, the $K$-module $\Lambda^2 \m_1$ is equivalent to the
adjoint representation of~$K_0$. Hence we have by Schur's Lemma that
\begin{equation*}
[\m_1, \m_1] \subseteq \k_0
\end{equation*}
and it follows that the space of $K$-equivariant maps $\Lambda^2
\m_1$ is real one-dimensional, since the adjoint representation
of~$\k_0$ is of real type.

We shall now explicitly describe the Lie bracket on the subalgebra
$\g_1 = \k \oplus \m_1$ of~$\g$. Let $e_1 , \ldots , e_n$ be an
orthonormal basis of~$\m_1$. Let $\Cl_n = \Cl(\m_1)$ be the {\em
Clifford algebra} over $\m_1$, i.e.\ the real algebra generated by
$e_1 , \ldots , e_n$ subject to the relation $e_i \cdot e_j + e_j
\cdot e_i = -2 \delta_{ij}$.

Identify $\k_0$ with the subset
\begin{equation}\label{SpinCliff}
\spin(n) = \spann_{\R} \left \{ e_i \cdot e_j \mid 1 \le i < j \le n
\right \} \subset \Cl_n,
\end{equation}
where the bracket is given by the commutator. Then the action
of~$\k_0$ on~$\m_1$ can be described by $[e_i \cdot e_j, e_k] = -e_k
\cdot e_i \cdot e_j$. Now let $\lambda \in \R$ be a constant. We
define a Lie bracket on $\m_1 \times \m_1$ by
\begin{equation}\label{CliffLie}
[e_i, e_j] = \lambda\, e_i \cdot e_j \in \k_0,\  i \neq j.
\end{equation}
It can be easily verified that the bracket defined in this way is a
skew symmetric map on $\g_1 \times \g_1 \to \g_1$ fulfilling the
Jacobi identity and indeed we have
\begin{equation*}
    \g_0 = \k_0 \oplus \m_1 \cong \left\{
                  \begin{array}{ll}
                    \spin(n+1) & \hbox{if $\lambda > 0;$} \\
                    \spin(n) \ltimes \R^n & \hbox{if $\lambda = 0;$} \\
                    \spin(n,1) & \hbox{if $\lambda < 0$.}
                  \end{array}
                \right.
\end{equation*}
Since the space of $K$-equivariant maps $\Lambda^2 \m_1 \to \k_0$ is
real one-dimensional and $K_1$ acts trivially on~$\m_1$, this
construction exhausts all possible Lie algebra structures on $\g_1 =
\k \oplus \m_1$. We have thus shown that the Lie bracket on $\m_1
\times \m_1$ is given by~(\ref{CliffLie}) for some real
constant~$\lambda$.

Let us now determine the Lie bracket on $\m_1 \times \m_2$. From
(\ref{CartanRel}) it follows that
$$
[\m_1, \m_2] \subseteq \m_2.
$$
Since we consider $K$-representations where the group $K$ acts on
$\m_1 \oplus \m_2$ as the automorphism group of a Clifford module of
the Clifford algebra~$\Cl_n$ (cf.\ \cite{riehm} or Propositions~3.2
and~3.3 of~\cite{ks}), we may assume that~$\m_2$ carries the
structure of a Clifford module over~$\Cl_n$ such that the action
of~$\k_0$ on~$\m_2$ is given by Clifford multiplication, i.e.\
$$
[X,w] = X \cdot w,\qquad \mbox{ for } X \in \spin(n),\, w \in \m_2,
$$
where we use the identification~(\ref{SpinCliff}). Let us determine
the space of $K$-equivariant maps $\m_1 \times \m_2 \to \m_2$.

We shall prove that $\Hom_K(\m_1\otimes\m_2,\m_2)$ has real
dimension at most one in the case of representations~\ref{r205} and
\ref{r202} and complex dimension at most one in case of
representations~\ref{r204} and \ref{r201}:

For the representations~\ref{r204}, where the action of~$zA$, $z \in
\C$, $|z|=1$, $A \in \Sp(n)$ on $\m_1 = \C$ is given by complex
multiplication with~$z^{2k}$, the real $K$-module $\m_1 \otimes
\m_2$ decomposes into two $4n$-dimensional submodules; on these two
submodules, the action of $zA$ is given by the matrices $z^{2k+1}A$
and $z^{2k-1}A$, respectively. Thus it follows that $\m_1 \otimes
\m_2$ contains a submodule equivalent to~$\m_2$ if and only if $k =
1$; furthermore, $\m_1 \otimes \m_2$ contains at most one such
submodule.

In case of the representations \ref{r205} and \ref{r202}, the real
$K$-representation $\m_1 \otimes \m_2$ contains at most one
submodule equivalent to~$\m_2$. (To verify this, one may for example
use Weyl's dimension formula to see that the tensor product of the
standard and the spin representation of $\Spin(3)$ or $\Spin(7)$
contains an irreducible module of real dimension $8$ or $48$,
respectively.) Since $\m_2$ is a $K$-representation of real type, it
follows that the space of $K$-equivariant maps $\m_1 \times \m_2 \to
\m_2$ is at most real one-dimensional.

Let us now consider the case of representation~\ref{r201}. Here the
real $K$-modules $\m_1$ and $\m_2$ are of real and complex type,
respectively, and hence we may view the $K$-representation $\m_1
\otimes \m_2$ as a complex $K$-module. The real tensor product $\m_1
\otimes \m_2$ decomposes into a complex 20-dimensional irreducible
plus a complex 4-dimensional module; thus the space of
$K$-equivariant maps $\m_1 \times \m_2 \to \m_2$ is at most complex
one-dimensional.

Clifford multiplication induces a map
\begin{equation*}
\m_1 \times \m_2 \to \m_2, \quad (v,w) \mapsto v \cdot w
\end{equation*}
which is $K$-equivariant since $K$ acts on $\m_1 \oplus \m_2$ as a
group of Clifford module automorphisms, cf. \cite{riehm}. We have
shown that the Lie bracket on $\m_1 \times \m_2$ is given by
\begin{equation}\label{CliffModLie1}
[v,w] = \mu \; v \cdot w
\end{equation}
where $\mu$ is a real constant in case of representations \ref{r205}
or \ref{r202} and a complex constant in case of
representations~\ref{r204} or \ref{r201}. Let $X,\,Y \in \m_1$, $Z
\in \m_2$. By using (\ref{CliffLie}) and (\ref{CliffModLie1}), it
follows from the Jacobi identity $ [[X,Y],Z] = [X,[Y,Z]] - [Y,[X,Z]]
$ that $\lambda = 2 \mu^2$.

We can always rescale the Lie bracket to have $\lambda=0$ or $\pm 1$
and furthermore change the sign of~$\mu$. In fact, consider a linear
map $f\colon\k\oplus\m_1\oplus\m_2\to\k\oplus\m_1\oplus\m_2$ mapping
$(k,x,v)\mapsto(k,\alpha x, v)$, and denote by $[\cdot,\cdot]'$ the
Lie bracket pulled back via $f$, i.e.\ $[x,y]'=f^{-1}[fx,fy]$. For
$y\in\m_2$ and $u,v\in\m_1$ we then have
\begin{equation}\label{rescliealg}\begin{split}
[u,v]' &=  f^{-1}[fu,fv]=\lambda fu\cdot fv = \alpha^2\lambda [u,v]' \in \k \quad \mbox{ and} \\
[u,y]' &=  f^{-1}[fu,fy]=f^{-1}(\mu fu\cdot fy) = \alpha\mu [u,y]' \in \m_2 .
\end{split}\end{equation}
We can thus rescale $(\lambda,\mu)$ to $(\alpha^2\lambda,\alpha\mu)$.
\medskip

We will now consider the Lie bracket on $\m_2 \times \m_2$. From
(\ref{CartanRel}) we know that
$$
[\m_2, \m_2] \subseteq \k \oplus \m_1.
$$

To proceed further, we distinguish between two cases, depending on
whether the modules $\m_1$ and $\m_2$ do or do not commute.


\paragraph{\em The case where $\m_1$ and $\m_2$ commute}
\label{CliffCommute}

Assume $[\m_1, \m_2] = 0$. Then we have $\lambda = 2 \mu^2 = 0$ and
it follows that $\m_1$ is an abelian ideal of~$\g$.
\par
If in addition $[\m_2, \m_2] = 0$ holds, then $\m = \m_1 \oplus
\m_2$ is an abelian ideal of~$\g$ and $M = G / K$ is a homogeneous
presentation of a Euclidean space.
\par
Now assume $0 \neq [\m_2, \m_2] \subseteq \m_1$; in this case, the
ideal $\m_1 \oplus \m_2$ of~$\g$ is a two-step nilpotent Lie algebra
with center $\m_1$. The Lie bracket on $\m_2 \times \m_2$ defines a
$K$-equivariant map $Z \mapsto J_Z$ from $\m_1$ to the skew
symmetric endomorphisms of~$\m_2$ by $\langle J_Z X \mid Y \rangle =
\langle Z \mid [ X, Y] \rangle$ and also a $K$-equivariant map
$\varphi \colon \m_1 \times \m_2 \to \m_2$ by
\begin{equation}\label{ind}
\langle \varphi(Z,X) \mid Y \rangle = \langle Z \mid [X,Y] \rangle =
\langle J_Z X \mid Y \rangle \quad \mbox{ for } Z \in \m_1,\;X,\,Y
\in \m_2.
\end{equation}
We have already shown that the space of $K$-equivariant maps $\m_1
\times \m_2 \to \m_2$ is real one-dimensional for representations
\ref{r205} or \ref{r202} and complex one-dimensional in case of
representations~\ref{r204} or \ref{r201}. One particular such
$K$-equivariant map is given by Clifford multiplication using the
inclusion $\m_1 \subset \Cl(\m_1)$. It follows that the map $J
\colon \m_1 \to \End (\m_2)$, $Z \mapsto J_{Z}$ is given by $J_Z
\colon X \mapsto \kappa\; Z \cdot X$ where $\kappa$ is a real
constant in case of representations \ref{r205} or \ref{r202} and a
complex constant in case of representations~\ref{r204} or
\ref{r201}, the dot denoting Clifford multiplication. (In case of
representations~\ref{r204} and \ref{r201}, $\m_2$ carries a complex
structure.) Since we assume $[\m_2,\m_2] \neq 0$, it follows that
$\kappa \neq 0$.

Now let $\rho = \sqrt{|\kappa|}$, $\varepsilon = \sgn(\kappa)$ in
case of representations \ref{r205} or \ref{r202} and let $\rho \in
\C$ be such that $\rho^2 = \kappa$ and $\varepsilon = +1$, in case
of representations~\ref{r204} or \ref{r201}. Define a bijective
linear map $f \colon \m \to \m$ by $f(Z) = \varepsilon Z$ for $Z \in
\m_1$; $f(X) = X / \rho$ for $X \in \m_1$. It follows that $\langle
f(Z) \mid [f(X),f(Y)] \rangle = \langle J_{f(Z)} f(X) \mid f(Y)
\rangle = \langle Z \cdot X \mid Y \rangle$. This shows that $\m$ is
isomorphic to the Lie algebra of a {\em generalized Heisenberg
group} of type
$$
  N(2,1),\; N(3;1,0),\; N(6,1),\; \mbox{or}\; N(7;1,0),
$$
cf.~\cite{btv}~3.1, and that $M$ endowed with the metric induced
from the scalar product $\langle \ \cdot \mid \cdot \ \rangle$ is
isometric to a generalized Heisenberg group. We have shown that all
possible two-step nilpotent Lie algebra structures on~$\m$ are given
by this construction.

Now assume $[\m_2, \m_2] \not \subseteq \m_1$. Since $\m_1$ is an
(abelian) ideal, we may consider the quotient algebra $\g / \m_1$,
of which $\k$ is a subalgebra via $\k \to \g / \m_1 \colon X \mapsto
X + \m_1$. From (\ref{CartanRel}) we have that $\k$ is a compact
symmetric subalgebra of~$\g / \m_1$. In case of representations
\ref{r204}, \ref{r201} and \ref{r202} this leads to a contradiction
with the classification of symmetric spaces~\cite{helgason}, since
in these cases the $K$-representation on~$\m_2$ is not equivalent to
an isotropy representation of any Riemannian symmetric space. For
representation~\ref{r205}, it follows that $\g / \m_1$ is a simple
Lie algebra of type~$\sp(n+1)$ or $\sp(n,1)$. But $\g / \m_1$ acts
nontrivially with nontrivial kernel $\k_1\oplus\m_2$ on $\m_1$ and
hence cannot be simple.


\paragraph{\em The case where $\m_1$ and $\m_2$ do not commute}
\label{CliffNotCommute}

In this case we have $\lambda = 2 \mu^2 \neq 0$. Furthermore, since
$\m_2$ is $K$-irreducible, we have $[\m_1, \m_2] = \m_2$
from~(\ref{CartanRel}). Let $G_1$ be the connected subgroup of~$G$
corresponding to~$\g_1$, which is closed by Theorem~\ref{g1}. Up to
scaling, the homogeneous space~$F = G_1 / K$ is (locally) isometric
to either $\HH^n$, in case $\lambda < 0$, or $\eS^n$, in case
$\lambda > 0$. By Theorem~\ref{nisomisotr} it follows that the
action of~$G_1$ on~$\m_2$ is isometric.

If $\lambda < 0$, it follows that $\k_0 \oplus \m_1 \cong
\spin(2,1)$, $\spin(3,1)$, $\spin(6,1)$, or $\spin(7,1)$; in
particular, $\k_0 \oplus \m_1$ is a simple real Lie algebra. It
follows that the corresponding group acts nontrivially and
isometrically on~$\m_2$, a contradiction. Thus we may henceforth
assume $\lambda > 0$. By the rescaling \eqref{rescliealg} we can
actually arrange for $\lambda= 1$ and $\mu = 1/\sqrt{2}$.

Now assume $[\m_2, \m_2] = 0$, i.e.\ $\m_2$ is an abelian ideal
of~$\g$. Then the Lie bracket of~$\g$ is uniquely determined
by~$\lambda$ and~$\mu$. Hence it suffices to exhibit examples of
such spaces $M$ in order to explicitly obtain the Lie algebra
structure (up to rescaling) and we conclude that $M = G / K$ is
isometric to one of the homogeneous spaces
$$
 \fracd{\Sp(1)\cdot\Sp(n)\ltimes\R^{4n}}{\U(1)\cdot\Sp(n)},\quad
 \fracd{\Sp(1)\cdot (\Sp(1) \cdot \Sp(n) \ltimes \R^{4n})}{\Delta
 \Sp(1) \cdot \Sp(n)},\quad
 \fracd{\Spin(7)\ltimes\delta_7}{\Spin(6)},\quad
 \fracd{\Spin(8)\ltimes\delta_8^+}{\Spin(7)}
$$
endowed with a $G$-invariant metric, where $\delta_7$ stands for the
$8$-dimensional spin representation of $\Spin(7)$ and $\delta_8^+$
denotes an $8$-dimensional half-spin representation of~$\Spin(8)$;
the~``$\Delta$'' indicates that the $\Sp(1)$-factor is diagonally
embedded into the two $\Sp(1)$-factors of~$G$.

If $0 \neq [\m_2, \m_2] \subseteq \k$, then $\k$ is symmetric in $\k
\oplus \m_2$. This is a contradiction for representations
\ref{r204}, \ref{r202} and \ref{r201}, since in these cases $\m_2$
is not equivalent to the isotropy representation of a Riemannian
symmetric space. For representation~\ref{r205}, it follows that
$[\m_2, \m_2] = \k$, but then we obtain from the Jacobi identity
\begin{equation}\label{Jaco}
[[\m_2, \m_2], \m_1] \subseteq [\m_2, [\m_2, \m_1]] = [\m_2, \m_2],
\end{equation}
the contradiction $\m_1 \subseteq \k$.

If $[\m_2, \m_2] = \m_1$, then (\ref{Jaco}) reduces to $[\m_1,\m_1]
\subseteq \m_1$, which implies $[\m_1, \m_1] = 0$. It follows that
$\m$ is a subalgebra of~$\g$ and since neither of the
$K$-irreducible spaces $\m_1$ or $\m_2$ is an ideal of~$\m$, the Lie
algebra~$\m$ must be simple. By \cite{helgason}, Corollary~II~6.5,
it follows that $\m$ contains a subalgebra isomorphic to~$\k$ whose
action on $\m$ is equivalent to the~$\k$-action. This is a
contradiction since $\m$ does not contain any $K$-submodule
equivalent to the adjoint representation of~$K$.

It only remains the case where $\m_1 \neq [\m_2, \m_2] \not
\subseteq \k$. In case $n = 6$,~$7$ it follows that $\g$ is simple,
immediately leading to a contradiction in case $n = 6$, since there
is no simple real Lie algebra of dimension~$29$; if $n = 7$ then
$\g$ is a simple real Lie algebra of dimension~$36$ containing
$\spin(7)$ as a subalgebra and it follows that $\g$ is isomorphic to
either $\spin(9)$, $\spin(8,1)$ or $\spin(7,2)$; however, only in
the first two cases the isotropy representation can be equivalent
to~\ref{r202}, hence
$$
  M = \Spin(9) / \Spin(7) \quad \mbox{or} \quad M = \Spin(8,1) /
  \Spin(7).
$$

Now consider the case of representations~\ref{r204} and \ref{r205}.
The subalgebras $\k_0 \oplus \m_1$ and $\k_1$ of~$\g_1$ commute and
hence are ideals of $\g_1$. Since we have already shown
that~$\lambda$ is positive, we know that $\k_0 \oplus \m_1$, and
therefore also $\g_1$, is a compact Lie algebra. By
(\ref{CartanRel}) we have that $\g_1$ is a compact symmetric
subalgebra of~$\g$.

In case of representation~\ref{r204} we have $\g_1 \cong \sp(1)
\oplus \sp(n)$ with $G_1$ acting transitively and almost effectively
on the unit sphere in~$\m_2$ and we obtain from the classification
of symmetric spaces~\cite{helgason} that $\g$ is isomorphic to
either $\sp(n +1)$ or $\sp(n,1)$. Hence it follows that
$$
 M = \Sp(n+1) / \U(1) \cdot \Sp(n) \quad \mbox{or} \quad
 M = \Sp(n,1) / \U(1) \cdot \Sp(n).
$$

In case of representation~\ref{r205}, we have that $\k_0 \oplus \m_1
\cong \so(4)$ and hence $\g_1 \cong \sp(1) \oplus \sp(1) \oplus
\sp(n)$. The compact symmetric subgroup $G_1$ acts transitively on
the unit sphere in~$\m_2$ and such that the action of the
$\Sp(n)$-factor is almost effective; it follows from the
classification of symmetric spaces~\cite{helgason} that $G / G_1$ is
a non-effective presentation of $\qH \HH^n$ or $\qH \P^n$. Hence
$\g$ is isomorphic to a Lie algebra direct sum $\sp(1) \oplus \sp(n
+ 1)$ or $\sp(1) \oplus \sp(n, 1)$. It follows that
$$
 M = \fracd{\Sp(1) \cdot \Sp(n+1)}{\Delta \Sp(1) \cdot \Sp(n)} \quad \mbox{or} \quad
 M = \fracd{\Sp(1) \cdot \Sp(n,1)}{\Delta \Sp(1) \cdot \Sp(n)},
$$
where the~``$\Delta$'' indicates that the $\Sp(1)$-factor is
diagonally embedded into the two simple factors of~$G$.
\end{proof}


    \subsection{A Lie algebra construction leading to
                homogeneous presentations of flat space}
    \label{aconstr}


In this section we will show that a Riemannian homogeneous space
whose isotropy representation is equivalent to~\ref{r203} from
Table~\ref{Tcoh2} is Euclidean.


\begin{proposition}\label{flat}
Let $M = G / K$ be a simply connected Riemannian homogeneous space
such that the action of~$K = \U(n)$ on the tangent space $\m = \m_1
\oplus \m_2 = \C \oplus \C^n$ of~$M$ is equivalent to
representation~\ref{r203}, i.e.\ $K$ acts by a power of the
determinant on~$\m_1$ and by the standard representation on~$\m_2$.
Then $M$ is isometric to Euclidean space.
\end{proposition}


\begin{proof}
It follows from Theorem~\ref{nisomisotr} that the action of~$G_1$
on~$\m_2$ is isometric, since the principal isotropy group of the
$K$-action on~$\m_2$ is isomorphic to~$\U(n-1)$. Assume the
$G_1$-action on~$\m_2$ is almost effective, then we have $G_1
\subseteq \SO(\m_2)$; however, since $\k = \u(n) \subset \so(2n)$ is
a maximal subalgebra, it follows that $\g_1 \cong \so(2n)$, which
leads to contradiction by a dimension count. Thus it follows that
the action of $\g_1$ has a two-dimensional kernel, which is an ideal
of~$\g_1$. Since the only two-dimensional $K$-invariant subspace
of~$\g_1$ is $\m_1$, this shows that $[\m_1, \m_2] = 0$.

Now observe that $\g_1 = \k \oplus \m_1$ is a subalgebra of~$\g$ by
Theorem~\ref{g1}. Since $\Lambda^2 \m_1$ is a one-dimensional
trivial $K$-module, it follows from Schur's Lemma that $ [\m_1,
\m_1] \subseteq \k_0, $ where $\k_0$ is the one-dimensional abelian
ideal of $\k =\u(n)$. But since we have already shown that $\m_1$ is
an ideal of~$\g_1$, it follows that $[\m_1,\m_1] = 0$.

Since the homogeneous space $\SO(2n) / \U(n) = \SO(\m_2) / K$ is
isotropy irreducible~\cite{wolf}, it is obvious by a dimension count
that the $K$-module $\m_1 \otimes \m_2$ does not contain any
submodule equivalent to~$\m_1$ or $\m_2$. Therefore $[\m_2, \m_2]
\subseteq \k$. Assume $[\m_2, \m_2] \neq 0$ and consider the Lie
algebra $\k \oplus \m_2$; it follows from the classification of
symmetric spaces~\cite{helgason} that $\k \oplus \m_2$ is a simple
real Lie algebra of type $\su(n+1)$ or $\su(n,1)$; however, $\k
\oplus \m_2$ acts on~$\m_1$ with nontrivial kernel~$\m_2$, a
contradiction.

We have proved that~$\m = \m_1 \oplus \m_2$ is an abelian ideal
of~$\g$ and hence $M$ is flat.
\end{proof}


    \subsection{Isotropy with a one-dimensional trivial module}
    \label{1triv1irr}


To complete the classification in the homogeneous case, it remains
to study manifolds $M = G / K$ where the isotropy representation has
a one-dimensional trivial submodule.

\begin{proposition}\label{1trivmod}
Let $M = G / K$ be a simply connected Riemannian homogeneous space
such that the action of~$K$ on the tangent space of~$M$ is of
cohomogeneity two and leaves a nonzero vector fixed. Assume that $M$
is not de~Rham decomposable and not symmetric. Then $M$ is isometric
to one of the following homogeneous spaces equipped with a
$G$-invariant metric:
$$
\frac{\SU(n+1)}{\SU(n)},\quad \frac{\SU(n,1)}{\SU(n)}, \quad N(1,k).
$$
\end{proposition}


\begin{proof}
Assume $\m_1$ is a one-dimensional trivial $K$-module and $\m_2$ is
an irreducible $K$-rep\-re\-sen\-ta\-tion of cohomogeneity one.
Clearly we have $[\m_1, \m_1] = 0$.


\paragraph{\em The case where $\m_1$ and $\m_2$ commute}

Assume $[\m_1, \m_2] = 0$. Then $\m_1$ lies in the center of~$\g$
and is hence an ideal of~$\g$. Consider the quotient algebra $\g_0 =
\g / \m_1$. It contains $\k$ as a subalgebra, acting irreducibly on
a $K$-invariant complement, denoted by~$\bar \m_2$, of~$\k$
in~$\g_0$; the action of $\k$ on~$\bar \m_2$ is clearly equivalent
to the $\k$-action on~$\m_2$, thus the corresponding homogeneous
space $M_0 = G_0 / K$ is strongly isotropy irreducible, with an
isotropy representation of cohomogeneity one. It follows that $M_0$
equipped with a $G_0$-invariant metric is either isometric to a
symmetric space of rank one or Euclidean.

Let $[\cdot, \cdot]_0$ be the Lie bracket of~$\g_0 = \k \oplus \bar
\m_2$. Consider the subspace $I = \bar \m_2 + [\bar \m_2, \bar
\m_2]_0$ of~$\g_0$. It follows from the Jacobi identity that $I$ is
an ideal of~$\g_0$, since
$$
[I, \k]_0 = [ \bar \m_2 + [\bar \m_2, \bar \m_2]_0, \k ]_0 \subseteq
[\bar \m_2, \k]_0 + [\bar \m_2, [\bar \m_2, \k]_0 ]_0 = I
$$
and
$$
[I, \bar \m_2]_0 \subseteq [\g_0, \bar \m_2]_0 = [\k, \bar \m_2]_0 +
[\bar \m_2, \bar \m_2]_0 = I.
$$
It follows that $\m_1 + \m_2 + [\m_2,\m_2]$ is an ideal of~$\g$.

Assume first the ideal~$I$ is simple. Then $\m_1 + \m_2 +
[\m_2,\m_2]$ is a central extension of the simple Lie algebra~$I$
and hence isomorphic to the Lie algebra direct sum of~$I$ and the
one-dimensional abelian Lie algebra~$\m_1$ by the Whitehead Lemma.
Since $\m_1$ lies in the center of~$\g$, it follows that also $\g$
is isomorphic to a Lie algebra direct sum of~$\g_0$ and a
one-dimensional abelian Lie algebra. Hence the homogeneous space $M
= G / K$ is de~Rham decomposable in this case.

Assume now $I = \bar \m_2$, i.e.\ $[\bar \m_2, \bar \m_2] \subseteq
\bar \m_2$. Then the $K$-irreducible space $I = \bar \m_2$ is either
an abelian or a simple ideal. We have already treated the case where
$I$ is simple, thus we may assume $I$ is abelian. In this case $\m =
\m_1 \oplus \m_2$ is a two-step nilpotent ideal of~$\g$ and it
follows that the Lie bracket on~$\m$ is given by a $K$-invariant
skew-symmetric bilinear form $\omega \colon \m_2 \times \m_2 \to \R$
such that
$$
[ X, Y ] = \omega( X, Y ) \, v\quad \mbox{ for all } X,Y \in \m_2,
$$
where $v \in \m_1$ is some nonzero vector. We may assume $\omega
\neq 0$ since otherwise $M$ is flat. The two-form $\omega$ defines a
nonzero skew-symmetric $K$-equivariant endomorphism $F$ on~$\m_2$ by
$\langle F (X) \mid Y \rangle = \omega ( X, Y )$. Hence $F^2$ is a
nonzero symmetric $K$-equivariant endomorphism on~$\m_2$ with at
least one negative eigenvalue. By Schur's Lemma it follows that $F^2
= -\lambda \, \id_{\m_2}$, for some $\lambda > 0$. Hence the
$K$-module $\m_2$ admits a complex structure defined by~$F /
\sqrt{\lambda}$. It follows that $M$ is isometric to a generalized
Heisenberg group of type
$$
  N(1,k),\;\mbox{where $k = \frac{1}{2} \dim_{\R} \m_2$.}
$$

We may now assume that $I \neq \bar \m_2$ and that $I$ is not
simple. Then $I$ contains a non-trivial proper ideal. Since $\k$
acts effectively on~$\bar \m_2$, the only possibility is that $\bar
\m_2$ is an ideal of~$I$, contradicting $I \neq \bar \m_2$.


\paragraph{\em The case where $\m_1$ and $\m_2$ do not commute}
In view of Theorem ~\ref{nisomisotr} we may assume that the
$G_1$-action on $\m_2$ is isometric. Furthermore, as $[\m_1, \m_2]
\neq 0$ we may assume that the group $G_1$ acts almost effectively,
i.e. with a discrete kernel $\Gamma\subset G_1$, on~$\m_2$. In fact,
if $\g_1$ had a kernel $\Lb$ on~$\m_2$ then, by the effectivity of
the $K$-action on~$\m_2$, $\Lb$ would be another $\Ad_K$-invariant
one-dimensional complement of~$\k$ in~$\g_1$; in this case we could
use $\Lb$ in place of $\m_1$ and would be in the previous case where
$\m_1$ and $\m_2$ commute.

Thus both $K$ and $G_1/\Gamma$ act effectively and transitively on
the sphere~$\eS\m_2$ in~$\m_2$. Therefore both must show up in the
classification of such groups. From Table~\ref{Tcoh1}, the only
cases where a group and a one-dimensional extension of it act
effectively and simply transitively on the same sphere are the cases
$$K= \SU(n),\quad G_1/\Gamma=\U(n),\quad \m_2=\C^n;$$
and
$$K= \Sp(n),\quad G_1/\Gamma=\U(1)\cdot\Sp(n),\quad \m_2=\qH^n. $$
In both cases, the only module isomorphic to
$\m_1\otimes\m_2\cong\m_2$ in $\g=\k\oplus\m_1\oplus\m_2$ is $\m_2$.
Therefore $G_1$ acts on $\m_2$ in the usual sense (i.e. the
projection \eqref{projg1act} is not needed). Also, $G_1$ contains an
element in its center acting as $-1$ on $\m_2$. Hence the Cartan
relations \eqref{CartanRel} hold and $(\g,\g_1)$ is an orthogonal
symmetric Lie algebra. It follows from the classification of
symmetric spaces that $G \cong \SU(n + 1)$ or $\cong \SU(n,1)$ if
$\m_2$ is not abelian, hence
$$
M = \SU(n + 1) / \SU(n) \quad \mbox{or} \quad M = \SU(n, 1) /
\SU(n).
$$
If $\m_2$ is abelian, then $\m=\m_1\oplus\m_2$ is an ideal of~$\g$.
Since the action of~$\m_1$ on~$\m_2$ commutes with the irreducible
$K$-action it follows that $\m_2$ carries a $K$-equivariant complex
structure and such that the action of elements in~$\m_1$ is given as
multiplication by purely imaginary elements of~$\C$. It follows as
in Theorem \ref{nisomisotr} that the simply connected Lie group
corresponding to $\m$ is isomorphic to the group $M = \R \ltimes
\C^n$, where the group multiplication is given by
$$
(t,v) \cdot (s,w) = (t+s, v +e^{it}w).
$$
Identifying $M$ with $\R^{1+2n}$, this induces a simply transitive
action of~$M$ on itself which is isometric if $M = \R^{1+2n}$
endowed with any Riemannian metric of the form
$$
 ds^2 = \lambda \, dx_1^2 + \mu \left ( dx_2^2 + \dots + dx_{1+2n}^2 \right ),
$$
where $\lambda$, $\mu > 0$ are constants. This shows that, for any
choice of a left invariant metric, $M = G / K$ is isometric to
Euclidean space, hence de~Rham decomposable.
\end{proof}


We have now completed the proof of Theorem~\ref{Coh2Result}.


    \section{The inhomogeneous case}
    \label{inhomclass}


In this section we classify simply connected manifolds admitting a
non-transitive group of isometries whose isotropy groups act with
cohomogeneity~$\leq 2$.


\begin{theorem}\label{InhomTh}
Let $M$ be a simply connected Riemannian manifold on which a closed
connected subgroup~$G$ of the isometry group acts non-transitively
and non-trivially. Let $Gp$, $p \in M$, be a principal orbit and
assume that the isotropy group $K = G_p$ acts with cohomogeneity two
on~$M$. Then the $G$-action has no exceptional orbits and at most
two singular orbits, which are fixed points. If there is a singular
orbit, then the principal orbits $G / K$ are spheres. Else the
principal orbits are simply connected rank-one symmetric spaces,
i.e.\ isometric to one of
\begin{equation}\label{rk1sym}
\R^n ,\ \eS^n,\ \C \P^n,\ \qH \P^n,\ \qO \P^2,\ \R \HH^n,\ \C
\HH^n,\ \qH \HH^n,\ \qO \HH^2.
\end{equation}
More specifically, see Table~\ref{Tcoh21}:
\begin{enumerate}
\item[(i)] If there are no singular orbits then $M$ is a warped
product $\R \times_f G/K$. \item[(ii)] If there is one singular
orbit then $M$ is isometric to $\R^{n+1}$ endowed with a
rotationally symmetric metric, i.e.\ a metric invariant under a
linear cohomogeneity one action of~$G$ on~$\R^{n+1}$. \item[(iii)]
If there are two singular orbits then $M$ is isometric to
$\eS^{n+1}\subset \R^{n+2}= \R^{n+1}\times \R$ carrying a metric
invariant under a linear cohomogeneity one action of $G$ on the
$\R^{n+1}$ factor.
\end{enumerate}
\begin{table}[h]
\begin{tabular}{|c|c|c|c|c|c|}
\hline No.\ & Singular orbits & $M$ &$G/K$  & Section $\Sigma$ &
Fundamental domain \\ \hline\hline (i) & none & $\R\times_f G/K$
&\eqref{rk1sym} & $\R$ & $\R$ \\ \hline (ii) & one point & $\approx
\R^{n+1}$    &$\eS^n$  & $\R$ & $\R^+_0$  \\ \hline (iii) & two
points & $\approx \eS^{n+1}$ &$\eS^n$ & $\eS^1$ & $[0,L]$ \\ \hline
\end{tabular}\vspace{1em}
\caption{Inhomogeneous manifolds with unit tangent bundle of
cohomogeneity one}\label{Tcoh21}
\end{table}
\end{theorem}
The metric on the warped product $\R\times_f G/K$ is $dt^2 \oplus
f(t)^2 g_0$ where $f$ is a smooth positive function on~$\R$ and
$g_0$ is a fixed $G$-invariant metric on~$G/K$. The metric $g_0$ is
uniquely determined up to a scaling factor since $G/K$ is isotropy
irreducible. All these manifolds are (degenerate) warped products,
i.e.\ admit an isometry $I\times_f G/K\to M$, where $I=\R, \R^+_0$
or some finite interval $[0,L]$, $L>0$, and $f>0$ on the interior of
$I$ and $f=0$ on $\partial I = \emptyset, \{0\}$ or $\{0,L\}$,
respectively.

\begin{proof}
We fix a point $p \in M$, such that $Gp$ is a principal orbit, and
denote by $K\subset G$ the isotropy subgroup at $p$. The isotropy
representation of~$K$ splits as
\begin{equation*}
T_p M = T_p Gp \oplus N_p Gp.
\end{equation*}
From this, it follows that $G$ acts with cohomogeneity one on~$M$
and also that $K$ acts with cohomogeneity one on~$T_p Gp$. Let
$\gamma \colon \R \to M$ be a normal geodesic such that
$\gamma(0)=p$ and $\dot{\gamma}(0) \perp T_{p} Gp$. Let $\Sigma =
\gamma(\R)$.

We prove the Theorem in the following steps:
\paragraph{{\bf 1.} {\em Each $G$-orbit intersects $\Sigma$ perpendicularly.}}
To prove this, let $Gq$ be a second $G$-orbit and choose a geodesic
$\psi$ realizing the distance of~$Gp$ to $Gq$, say $\psi(0)=gp$ for
some $g \in G$. This geodesic meets both $Gp$ and $Gq$
perpendicularly. Thus $g^{-1}\psi\subset \Sigma$ and $\Sigma$ meets
$Gp$ and $Gq$ perpendicularly.

\paragraph{{\bf 2.}
{\em There are no exceptional orbits. If there are singular orbits,
then these are fixed points and the principal orbits are spheres.
Else there are no singular orbits and the principal orbits are
simply connected.}} Let $G_q\subset G$ be the isotropy group at
$q\in\Sigma$. Hence $K \subseteq G_q$. Since $\g = \k \oplus T_pGp$
and $T_pGp$ is
irreducible, the Lie algebra $\g_q$ of $G_q$ is either $\g$ or $\k$.\\
a) If $\g_q=\g$ then $G_q=G$ and $q$ is a fixed point.         \\
b) If $\g_q=\k$ then $K$ has finite index $\iota$ in $G_q$. Assume
$Gq$ is an exceptional orbit i.e.\ $\iota > 1$. The map $Gp=G/K\to
Gq=G/G_q$, $gp\mapsto gq$, is well defined and has fibre $G_q/K$
with cardinality $\#G_q/K=\iota $. Thus it is a covering with
$\iota$ sheets. Since $Gp$ is connected, $\pi_1(Gp)$ is a subgroup
of positive index $\iota$ in $\pi_1(Gq)$. We may assume $p$
arbitrarily close to $q$ since the union of the principal orbits is
dense. By the Slice Theorem we may identify $Gp$ with a submanifold
of the sphere normal bundle $\eS \nu(Gq)$ of~$Gq$. Since $Gp$ is
connected we have
$$
Gp = G/K= \eS\nu(Gq) \to Gq = G/G_q
$$
a two-fold covering. We now apply the Theorem of Seifert-van Kampen
to the decomposition
$$
M = \overline{D\nu(Gq)} \cup M\setminus D\nu(Gq)
$$
where $D\nu(Gq)$ is a $G$-invariant open tubular neighborhood of
$Gq$ from the Slice Theorem. The intersection is
$$
\overline{D\nu(Gq)}\cap M\setminus D\nu(Gq) = \eS\nu(Gq) = Gp.
$$
Thus $\pi_1(M)\cong\pi_1(Gq) * \pi_1(M\setminus D\nu(Gq))
/\pi_1(Gp)$ cannot be trivial because
$\pi_1(Gp)\to\pi_1(Gq)\cong\pi_1(\overline{D\nu(Gq)}) $ is not
surjective. Thus there are no exceptional orbits and singular orbits
are fixed points.

If there is a fixed point then the distance sphere of a fixed point
is a principal orbit. If there is no fixed point then all orbits are
principal. Thus $M$ is a $G/K$-bundle over a one-dimensional
manifold which is either $\R$ or $\eS^1$. If it were $\eS^1$ or if
the fibres $G/K$ were not simply connected then so would be $M$, by
the homotopy sequence $0=\pi_2(G\backslash
M)\to\pi_1(G/K)\to\pi_1(M) \to\pi_1(G\backslash M)\to\pi_0(G/K)=0 $
of the fibre bundle $G/K \to M \to G\backslash M $.


\paragraph{{\bf 3.} {\em $M$ is a (possibly degenerate) warped product.}}
For $g \in G$ we have that $gp \in \Sigma$ implies $g\Sigma=\Sigma$,
because $\Sigma$ is the image of a geodesic $\gamma_v$ emanating
from $p=\gamma_v(0)$ with $\dot{\gamma}_v(0)=v \in
T_pGp^\perp=T_p\Sigma $. Hence $gv \in T_qGp^\perp=T_q\Sigma$ and
therefore $gv=\pm\dot{\gamma}_v(t_q)$ if $q=\gamma_v(t_q)$. Let
$W=\{g \in G\mid g\Sigma=\Sigma\}$.

If there are no fixed points then $M$ is a fibre bundle over
$G\backslash M=W\backslash \Sigma$. Since $M$ is simply connected we
have $W\backslash \Sigma=\Sigma=\R$ and the $G$-orbits intersect
$\Sigma$ only once. We thus have a diffeomorphism
\begin{equation}\label{isometrwp}
\begin{split}
\Sigma \times Gp & \to  M, \\
(s,gp) &\mapsto  gs.
\end{split}
\end{equation}
Since the $G$-invariant metric $g_0$ on~$G/K$ is unique up to
scaling we can find a positive function~$f$ on~$\Sigma$ such that
\eqref{isometrwp} becomes an isometry when $\R\times G/K=\Sigma
\times Gp$ is endowed with the warped product $dt^2\oplus f(t)^2
g_0$. Thus $M$ is of the first type in Table~\ref{Tcoh21}.

If there is a fixed point $q$ there is $w \in W\setminus K$ with
$w^2 \in K$ and $W$ is generated by $w$ and $K$. If $\Sigma=\R$
there are no other fixed point and if $\Sigma=\eS^1$ we have at most
two fixed points. In both cases the principal orbits are spheres.
Away from the fixed points $M$ is a warped product over the interior
of a fundamental domain $\Sigma_0 \subset \Sigma$ for the $W$-action
on~$\Sigma$. This can only be $\Sigma_0=\R^+$ or $\Sigma_0=(0,L)$
for some $L>0$, which produces the second and the third type in
Table~\ref{Tcoh21}, respectively. One has to replace $\R= \Sigma$ in
\eqref{isometrwp} by the fundamental domain $\Sigma_0=\R^+$ or
$\Sigma_0=(0,L)$.
\end{proof}


\end{document}